\documentclass[11pt]{article}
\usepackage[greek,english]{babel}
\usepackage[iso-8859-7]{inputenc}
\usepackage{amssymb}
\usepackage{amsmath}
\usepackage{amsthm}
\usepackage{latexsym}
\usepackage{amsfonts}
\usepackage{graphicx}
\usepackage{graphics}

\theoremstyle{plain}
\newtheorem{thm}{Theorem}[section]   
\newtheorem{prop}[thm]{Proposition}
\newtheorem{lem}[thm]{Lemma}
\newtheorem{cor}[thm]{Corollary}

\theoremstyle{definition}

\newtheorem*{Proof}{Proof}

\newcommand{\ga} {{\gamma}}
\newcommand{\al} {{\alpha}}
\newcommand{\Ga} {{\varGamma}}

\newcommand{\ld} {{\ldots}}
\newcommand{\sm} {{\smallsetminus}}
\newcommand{\thi} {{\theta}}

\newcommand{\de} {{\delta}}

\newcommand{\la} {{\lambda}}

\newcommand{\vPsi} {{\varPsi}}

\newcommand{\Fi} {{\varPhi}}
\newcommand{\e} {{\varepsilon}}
\newcommand{\f} {{\varphi}}

\newcommand{\dis}{\displaystyle}

\newcommand{\bfi}{\bar{f}}
\newcommand{\ct}{{\cal{T}}}

\newcommand{\cu}{{\cal{U}}}

\newcommand{\ch}{{\cal{H}}}

\newcommand{\cc}{{\cal{C}}}

\newcommand{\tf}{{\widetilde{f}}}

\newcommand{\tg}{{\widetilde{G}}}
\newcommand{\tc}{{\widetilde{C}}}
\newcommand{\mfd}{{\mathfrak{D}}}
\newcommand{\ra}{{\rightarrow}}
\newcommand{\oD}{{\overline{D}}}
\newcommand{\oG}{{\overline{G}}}

\newcommand{\para}{{\parallel}}

\newcommand{\qb}{$\quad\blacksquare$}
\def\1{\it1\hspace*{-0.150cm}{\footnotesize{I}}}

\def\R{{\mathbb{R}}}
\def\C{{\mathbb{C}}}

\def\N{{\mathbb{N}}}

\begin{document}
\title{\bf Uniform convergence of translation operators}
\author{\bf N. Tsirivas}
%
\date{}
\maketitle
\noindent
\begin{abstract}
Let $\thi$ be a fixed positive number, $\thi\in(0,1)$ and let $\la=(\la_n)_{n\in\N}$ be a fixed sequence of non-zero complex numbers, so that $\la_n\ra\infty$. We shall apply the functions $g_n:[0,\thi]\times\C\ra\C$, defined as $g_n((t,z))=z+\la_ne^{2\pi it}$ for  each $(t,z)\in[0,\thi]\times\C$.

We shall consider the space $C([0,\thi]\times\C)$ of continuous functions on $[0,\thi]\times\C$, as endowed with the topology of uniform convergence on compacta and let $\rho$ be the usual metric in $C([0,\thi]\times\C)$. For an entire function $f\in\ch(\C)$ we shall denote that
\[
\bfi:[0,\thi]\times\C\ra\C, \ \ \bfi((t,z))=f(z) \ \ \text{for every} \ \ (t,z)\in[0,\thi]\times\C.
\]
We will prove that the equation: $\dis\lim_{n\ra+\infty}\rho((x\circ g_{y_n},\bfi))=0$ does not have any solution $(x,y_n)$ where $x\in\ch(\C)$ and $y_n$ is an strictly increasing subsequence of natural numbers and $f\in\ch(\C)$ is a given non-constant entire function. When $f$ is a constant entire function, then the above equation has infinitely several solutions, according to a result provided by G. Costakis.
\end{abstract}%
{\em Keywords}\,: hypercyclic operator, common hypercyclic vectors, translation operator. \smallskip\\
MSC (2020) 47A16

\section{Introduction}\label{sec1}
\noindent

A classical result of Birkhoff \cite{2}, which goes back to 1929, says that there are entire functions of which the integer translates are dense in the space of the entire functions endowed with the $\ct_u$ topology of local uniform convergence (see also Luh \cite{8} for a more general statement). Birkhoff's proof was constructive.

Much later, during the 80's, Gethner and Shapiro \cite{6} and independently Grosse-Erdmann \cite{7} showed that Birkhoff's result can be recovered as a particular case of a much more general theorem, through the use of Baire's category theorem.

This approach simplified Birkhoff's argument substantially and, in addition, it gave us precise information on the topological size of these functions. In particular, Grosse-Erdmann proved that for every fixed sequence of complex numbers $(w_n)$ with $w_n\ra\infty$, the set
\[
\big\{f\in\ch(\C)|\overline{\{f(z+w_n):n\in\N\}}=\ch(\C)\big\}
\]
is $G_\de$ and dense in $\ch(\C)$, and hence ``large'' in terms of topology.

More recently, Costakis and Sambarino \cite{4} established a notable strengthening of Birkhoff's result. Namely, they showed that, for almost all entire functions $f$, in the context of Baire's category theorem, the set of the translates of $f$ with respect to $na$, $n\in\N$, is dense in the space of all entire functions for every non-zero complex number $\al$. The significant new element here is the uncountable range of $\al$.

More specifically, the set
\[
\bigcap_{\al\in\C\sm\{0\}}\big\{f\in\ch(\C)|\overline{\{f(z+n\al):n\in\N\}}=\ch(\C)\big\}
\]
is residual in $\ch(\C)$, that is, it contains a $G_\de$ and dense subset of $\ch(\C)$.

In particular, it is non-empty. Subsequently, Costakis \cite{3} examined a similar result, where $n$ can be replaced by more general sequences $(\la_n)$ of non-zero complex numbers, so that $\la_n\ra\infty$.

In this direction, Costakis \cite{3} proved that if the sequence $(\la_n)$ satisfies a certain condition, then the desired conclusion is reached, provided that we focus on $\al\in C(0,1):=\{z\in\C/|z|=1\}$. Under this condition, it is proved that the set
\[
\bigcap_{\al\in C(0,1)}\big\{f\in\ch(\C)|\overline{\{f(z+\la_n\al):n\in\N\}}=\ch(\C)\big\}
\]
is residual in $\ch(\C)$.

The proof of this result follows a similar method to the one used to prove a similar result in \cite{3}. In the same article \cite{3}, Costakis examined a simpler and more specific case of the above result. In particular:\\
In the above set, the request is to find an entire function $f$, so that:
\[
\overline{\{f(z+\la_n\al):n\in\N\}}=\ch(\C) \ \ \text{for every} \ \ \al\in C(0,1).
\]
Let $C\subseteq\ch(\C)$ be the set of constant functions of $H(\C)$, that is
\[
C=\big\{g:\C\ra\C|\ \ \text{there exists} \ \ \al\in\C \ \ \text{such that} \ \ g(z)=\al \ \ \text{for every} \ \ z\in\C\big\}
\]
Firstly, Costakis proved in \cite{3} that there is some $f\in\ch(\C)$, so that $C$ is a subset of $\overline{\{f(z+\la_n\al):n\in\N\}}$ for every $\al\in C(0,1)$.

This result offered a different and notable proof. In fact, he proved something stronger in this case, by adding a stronger condition of convergence. More specifically, Costakis proved the following result:

Let $(\la_n)_{n\in\N}$ be a sequence of non-zero complex numbers, so that $\la_n\ra\infty$. We have the set:\\
$\cu_\C(\la)=\big\{f\in\ch(\C)|$ for every $x,y\in\R$ such that $0<x<y$, for every $\thi\in(0,1)$ and for every $\al\in\C$ there is a sequence $m=(m_n)_{n\in\N}$, so that $m_n\in\{\la_n,n\in\N\}$, for every $n\in\N$, so that, for every compact subset $L\subseteq\C$ $\dis\sup_{r\in[x,y]}\dis\sup_{t\in[0,\thi]}\dis\sup_{z\in L}|f(z+m_nre^{2\pi it})-\al|\ra0$ as $n\ra+\infty\big\}$.

Costakis \cite{3} proved that the above set $\cu_\C(\la)$ is a $G_\de$ and dense subset of $\ch(\C)$. However, he did not use this method in the general case and gave a completely different proof in the general case.

Therefore, it is reasonable to ask if we can deal with the general case by imitating the proof of the above specific case. In this paper, we shall prove that this cannot be done. More specifically, we will prove here the following result:

Let $\la=(\la_n)_{n\in\N}$ be a given sequence of non-zero complex numbers, so that $\la_n\ra\infty$, $\thi\in(0,1)$ and $G\in\ch(\C)$, where $\thi$ is a given number and $G$ is a given non constant function. We shall consider the set $\cu(\la,\thi,G)=\big\{f\in\ch(\C)|$ there exists a sequence $m=(m_n)_{n\in\N}$, where $m_n\in\{\la_n,n\in\N\}$ for every $n\in\N$, so that for every compact set $L\subseteq\C$ $\dis\sup_{(t,z)\in[0,\thi]\times L}|f(z+m_ne^{2\pi it})-G(z)|\ra0$ as $n\ra+\infty\big\}$.

Our main result is that $\cu(\la,\thi,G)=\emptyset$, that confirms that we cannot achieve the general result of Costakis \cite{3} by giving a proof similar to the proof of the specific case of constant functions. Let us present now our result here in a more formal setting. We fix a positive number $\thi\in(0,1)$ and a sequence $\la=(\la_n)_{n\in\N}$ from non-zero complex numbers, so that $\la_n\ra\infty$. For every $n\in\N$ we have the function $h_n:[0,\thi]\times\C\ra\C$, so that $h_n((t,z))=z+\la_ne^{2\pi it}$ for every $(t,z)\in[0,\thi]\times\C$. For every $n\in\N$ we use the map $T_n:\ch(\C)\ra C([0,\thi]\times\C)$, so that $T_n(f)=f\circ h_n$, for every $n\in\N$, $f\in\ch(\C)$. The sequence $(T_n)_{n\in\N}$ is a sequence of linear and continuous operators.

We consider the spaces $\ch(\C)$ and $C([0,\thi]\times\C)$ as endowed with the topology of uniform convergence on compacta. For every $f\in\ch(\C)$ we assume the orbit of $f$ under the sequence $T_n$, $n\in\N$ to be the set
$O(f,T_n)=\{g\in C([0,\thi]\times\C)|$ there is $n\in\N$, so that $g=T_n(f)\}$. For every $f\in\ch(\C)$ we use the function $\tf:[0,\thi]\times\C\ra\C$ so that $\tf((t,z))=f(z)$ for every $(t,z)\in[0,\thi]\times\C$. Our main result suggest that if $G\in\ch(\C)\sm C$, then $\tg\notin O(f,T_n)'$ for every $f\in\ch(\C)$, where with $O(f,T_n)'$ we refer to the set of accumulation points of $O(f,T_n)$ and with $C$ we refer to the set of constant functions of $\ch(\C)$.

Let
\[
\tc=\big\{g\in C([0,\thi]\times\C)| \;\;\text{there is a} \;\; f\in C\;\; \text{so that} \;\; g=\tf\big\} \ \ \text{and}
\]
\[
R=\big\{g\in C([0,\thi]\times\C)| \;\;\text{there is a} \;\; f\in\ch(\C),\;\; \text{so that} \;\; g=\tf\big\}.
\]
Based on our main result in Theorem \ref{thm4.2}, Lemma \ref{lem3.1} and the result  of Costakis \cite{3} we shall have:
\[
\bigg(\bigcup_{f\in\ch(\C)}O(f,T_n)\bigg)\cap R=\tc  \ \ \text{and}  \ \ \bigg(\bigcup_{f\in\ch(\C)}(O(f,T_n)')\bigg)\cap R=\tc \ \ \text{and}
\]
\[
O(f,T_n)\cap\big(O(f,T_n)'\big)=\emptyset \ \ \text{for every} \ \ f\in\ch(\C)\sm C.
\]

The paper is organized as follows: \\
After the introduction in Section \ref{sec1}, we shall prove Proposition \ref{prop2.1} that is a specific case of our main result, in the case that $G\in\ch(\C)$ is not a constant function, so that $G(0)=0$ and $G'(0)\neq0$.

In order to prove Proposition \ref{prop2.1}, we shall use 5 lemmas.

In Section \ref{sec2}, we shall analyze the proofs of the 5 lemmas.

In Section \ref{sec3}, we shall give a helping corollary and the proof of our main result in Theorem \ref{thm4.2}.

There are several results concerning the existence or non-existence of common hypercyclic vectors for translation operators, see \cite{10}, \cite{11}, \cite{1}, \cite{9}, \cite{5}.
\section{A specific case}\label{sec2}
\noindent

We shall use the $\N$ and $\C$ abbreviations for these sets of natural and complex numbers respectively.

We fix a positive number $\thi$ and a sequence $\la=(\la_n)_{n\in\N}$ from complex numbers.\\
We define $\ch(\C)$ for the set of entire functions. We fix $G\in\ch(\C)$.

Let's consider the set:\\
$\cu(\la,\thi,G)=\big\{f\in\ch(\C)|$ there exists a sequence $m=(m_n)_{n\in\N}$, where $m_n\in\{\la_n,n\in\N\big\}$ for every $n\in\N$, such that for every compact set $L\subseteq\C$\\
$\dis\sup_{(t,z)\in[0,\thi]\times L}|f(z+m_ne^{2\pi it})-G(z)|\ra0$ as $n\ra+\infty\big\}$.

We use some notations for reasons of simplification:\\
We set: $K=\big\{z\in\C|$ there is $\al\in[0,\thi]$ so that $z=e^{2\pi\al i}\big\}$.\\
We fix a sequence $m=(m_n)_{n\in\N}$ from complex numbers.\\
We consider the functions $\f_n:K\times\C\ra\C$, so that $\f_n((t,z))=z+m_nt$ for every $(t,z)\in K\times\C$, $n=1,2,\ld$, and their partial functions $\f^t_n:\C\ra\C$, so that $\f^t_n(z)_=\f_n((t,z))$ for every $n=1,2,\ld$, $t\in K$, $z\in\C$.\\
The functions $\f^t_n$ are obviously entire for every $n\in\N$, $t\in K$.

We assume $f\in\ch(\C)$.\\
Let's now consider the functions: $\Fi_n(f):K\times\C\ra\C$, where $\Fi_n(f)=f\circ\f_n$ and $\Fi^t_n(t)=f\circ\f^t_n$, for every $n\in\N$, $t\in K$ and $f\in\ch(\C)$. The functions $\Fi^t_n(f)$ are obviously entire for every $t\in K$, $n\in\N$, $f\in\ch(\C)$.

Afterwards, we consider the functions: $\vPsi_n(f):K\times\C\ra\C$, where
\[
\vPsi_n(f)((t,z))=\int_{[0,z]}(\Fi^t_n(f))(j)dj, \ \ \text{for every} \ \ n\in\N, \ \ f\in\ch(\C), \ \ (t,z)\in K\times\C.
\]
If $f\in\ch(\C)$, we set $f_0$ for the entire function, so that:
\[
f_0(z)=\int_{[0,z]}f(j)dj, \ \ \text{for every} \ \ z\in\C.
\]
The function $f_0$ is the unique anti-derivative of $f$, so that $f_0(0)=0$. After setting the above notations and symbolisms, we proceed to the following Proposition \ref{prop2.1}.\
\begin{prop}\label{prop2.1}
Let $\la=(\la_n)_{n\in\N}$ be a sequence of non-zero complex numbers, so that $\la_n\ra\infty$. We assume also that $\thi\in(0,1)$ and $G\in\ch(\C)$, so that $G(0)=0$ and $G'(0)\neq0$. Then we have: $\cu(\la,\thi,G)=\emptyset$.
\end{prop}
\begin{Proof}
So, as to provide a proof by contradiction, we suppose that $\cu(\la,\thi,G)\neq\emptyset$. Let $F\in\cu(\la,\thi,G)$. Then by Lemma \ref{lem3.1} there is a subsequence $m=(m_n)_{n\in\N}$ of $\la$, from different terms so that for every compact subset $L$ of $\C$
\begin{eqnarray}
\sup_{(t,z)\in[0,\thi]\times L}|F(z+m_ne^{2\pi it})-G(z)|\ra0 \ \ \text{as} \ \ n\ra+\infty. \label{prop1}
\end{eqnarray}
Based on Lemma \ref{lem3.2}, we assume that for every compact subset $L$ of $\C$
\begin{eqnarray}
\sup_{(t,z)\in[0,\thi]\times L}|F'(z+m_ne^{2\pi it})-G'(z)|\ra0 \ \ \text{as} \ \ n\ra+\infty. \label{prop2}
\end{eqnarray}
We set $f=F'$.

Of course, we have $G(z)=\dis\int_{[0,z]}G'(j)dj$ for every $z\in\C$, because $G(0)=0$. According to Lemma \ref{lem3.4}, we have that for every compact subset $L$ of $\C$
\begin{eqnarray}
\sup_{(t,z)\in K\times L}|\vPsi_n(f)((t,z))-F(z)|\ra0 \ \ \text{as} \ \ n\ra+\infty.  \label{prop3}
\end{eqnarray}
We assume the partial functions $\vPsi_n(f)^t:\C\ra\C$, for $n=1,2,\ld$, so that:
\begin{eqnarray}
\vPsi_n(f)^t(z)=\vPsi_n(f)((t,z)), \ \ \text{for every} \ \ n=1,2,\ld, \ \ t\in K \ \ \text{and} \ \ z\in\C.  \label{prop4}
\end{eqnarray}
We fix $t_0\in K$ and $n_0\in\N$. Then, by (\ref{prop4}) we have:
\begin{eqnarray}
(\vPsi_{n_0}(f)^{t_0})'(z)=f(z+m_{n_0}t_0) \label{prop5}
\end{eqnarray}
(as it is well-known by complex analysis and the definition of function $(\vPsi_{n_0}(f)^{t_0})$).

According to the chain rule we have
\begin{eqnarray}
(F\circ\f_{n_0}^{t_0})'(z)=f(z+m_{n_0}t_0) \ \ \text{for every} \ \ z\in\C.  \label{prop6}
\end{eqnarray}
Based on equalities (\ref{prop5}), (\ref{prop6}) we take that there is a constant $\cc_{n_0,t_0}\in\C$, that is dependent on $n_0$ and $t_0$, so that
\begin{eqnarray}
\vPsi_{n_0}(f)^{t_0}(z)=(F\circ\f^{t_0}_{n_0})(z)+\cc_{n_0,t_0} \ \ \text{for every} \ \ z\in\C.  \label{prop7}
\end{eqnarray}
Based on equality (\ref{prop7}) we have:
\begin{align}
&\vPsi_{n_0}(f)^{t_0}(0)=(F\circ\f^{t_0}_{n_0})(0)+\cc_{n_0,t_0}\;\;\text{that implies}\;\; F(m_{n_0}t_0)+\cc_{n_0,t_0}=0,\nonumber \\
&\text{so} \;\;\cc_{n_0,t_0}=-F(m_{n_0}t_0).  \label{prop8}
\end{align}
We have $F'=f'_0=f$. So, there is a constant $c\in F$, so that: $F(z)=f_0(z)+c$ for every $z\in\C$. Thus, we have $c=F(0)$ and $F(z)=f_0(z)+F(0)$ for every $z\in\C$, or
\begin{eqnarray}
F(z)=F(0)+\int_{[0,z]}f(j)dj \ \ \text{for every} \ \ z\in\C.  \label{prop9}
\end{eqnarray}
Based on equalities (\ref{prop8}) and (\ref{prop9}), we have
\begin{eqnarray}
\cc_{n_0},t_0=-\bigg(F(0)+\int_{[0,m_{n_0}t_0]}f(j)dj\bigg).  \label{prop10}
\end{eqnarray}
Based on equalities (\ref{prop7}) and (\ref{prop10}) we have:
\begin{eqnarray}
\vPsi_{n_0}(f)^{t_0}(z)=(F\circ\f^{t_0}_{n_0})(z)-F(0)-\int_{[0,m_{n_0}t_0]}f(j)dj \ \ \text{for every} \ \ z\in\C.  \label{prop11}
\end{eqnarray}
Based on (\ref{prop11}) we have that for every $t\in K$ for every $n\in\N$, $z\in\C$
\begin{eqnarray}
\vPsi(f)((t,z))=F(z+m_nt)-F(0)-\int_{[0,m_nt]}f(j)dj.  \label{prop12}
\end{eqnarray}
Based on (\ref{prop3}) and (\ref{prop12}) we have: for every compact subset $L$ of $\C$
\begin{eqnarray}
\sup_{(t,z)\in}\sup_{K\times L}|F(z+m_nt)-F(0)-\int_{[0,m_nt]}f(j)dj-G(z)| \ \ \text{as} \ \ n\ra+\infty.  \label{prop13}
\end{eqnarray}
By using the triangle inequality, we take that for every $n\in\N$, $L\subseteq\C$, $L$ compact
\begin{align}
\sup_{t\in K}\bigg|\int_{[0,m_nt]}f(z)dz+F(0)\bigg|\le&\sup_{(t,z)\in K\times L}\bigg|F(z+m_nt)-F(0)-\int_{[0,m_nt]}f(j)dj-G(z)\bigg| \nonumber\\
&+\sup_{(t,z)\in K\times L}|F(z+m_nt)-G(z)|.
\label{prop14}
\end{align}
By convergence of (\ref{prop1}) and (\ref{prop13}) based on the inequality (\ref{prop14}), we have
\begin{eqnarray}
\sup_{t\in K}\bigg|\int_{[0,m_nt]}f(z)dz+F(0)\bigg|\ra0 \ \ \text{as} \ \ n\ra+\infty.  \label{prop15}
\end{eqnarray}
By convergence of (\ref{prop15}) we also have:
\begin{eqnarray}
\bigg|\int_{[0,m_n]}f(z)dz+F(0)\bigg|\ra0 \ \ \text{as} \ \ n\ra+\infty.  \label{prop16}
\end{eqnarray}
Based on the properties of complex integral and triangle inequality we have:
\begin{align}
\sup_{t\in K}\bigg|\int_{[m_n,m_nt]}f(z)dz\bigg|\le&\sup_{t\in K}\bigg|\int_{[0,m_nt]}f(z)dz+F(0)\bigg| \nonumber \\
&+\bigg|\int_{[0,m_n]}f(z)dz+F(0)\bigg|, \ \ \text{for every} \ \ n\in\N.  \label{prop17}
\end{align}
By (\ref{prop15}), (\ref{prop16}) and (\ref{prop17}) we have
\begin{eqnarray}
\sup_{t\in K}\bigg|\int_{[m_n,m_nt]}f(z)dz\bigg|\ra0 \ \ \text{as} \ \ n\ra+\infty. \label{prop18}
\end{eqnarray}
By (\ref{prop18}) we have
\begin{eqnarray}
\int_{[m_n,m_ne^{2\pi i\thi}]}f(z)dz\ra0 \ \ \text{as} \ \ n\ra+\infty.  \label{prop19}
\end{eqnarray}
By (\ref{prop2}) we have: for $L=\{0\}$
\begin{eqnarray}
\sup_{t\in[0,\thi]}|f(m_ne^{2\pi it})-G'(0)|\ra0 \ \ \text{as} \ \ n\ra+\infty.  \label{prop20}
\end{eqnarray}
By (\ref{prop19}), (\ref{prop20}) and Lemma \ref{lem3.5} we reach to a contradiction and the proof of Proposition \ref{prop2.1} is complete.\qb
\end{Proof}

In the following pages, we shall prove the lemmas we have used in the above Proposition \ref{prop2.1}.
\section{Proofs of 5 Lemmas}\label{sec3}
\begin{lem}\label{lem3.1}
Let $\la=(\la_n)_{n\in\N}$ be a sequence of non-zero complex number and $\thi$ be a positive number. We suppose that $f\in\cu(\la,\thi,G)$, where $\cu(\la,\thi,G)$ is the set that is defined in Proposition \ref{prop2.1}. Then, the sequence $m=(m_n)_{n\in\N}$, which satisfies the condition of $\cu(\la,\thi,G)$, that is for every compact subset $L\subseteq G$,\\
$\dis\sup_{(t,z)\in[0,\thi]\times L}|f(z+m_ne^{2\pi it})-G(z)|\ra0$ as $n\ra+\infty$, is an infinite subset of $\C$ and can be chosen to be a subsequence of $\la$ from different terms.
\end{lem}
\begin{Proof}
We set
\[
a_n=\sup_{(t,z)\in[0,\thi]\times L}|f(z+m_ne^{2\pi it})-G(z)|, \ \ n=0,1,2,\ld,
\]
for some specific compact subset $L\subseteq\C$.

We suppose that $a_n=0$ for some $n\in\N$. Let $z_0\in L$. Because of $a_n=0$, we have $f(z_0+m_ne^{2\pi it})=G(z_0)$ for every $t\in[0,\thi]$. Because of $m_n\neq0$, based on our hypothesis, we take that the function $f$ is a constant by the principle of analytical continuation, so we have $f(z)=G(z_0)$ for every $z\in\C$. As a result,
\[
|f(z+m_ne^{2\pi it})-G(z)|=|G(z)-G(z_0)| \ \ \text{for every} \ \ z\in\C, \ \ n\in\N \ \ t\in[0,\thi] \ \ \text{and}
\]
\[
a_n=\sup_{(t,z)\in[0,\thi]\times L}|f(z+m_ne^{2\pi it})-G(z)|=\sup_{z\in L}|G(z)-G(z_0)|
\]
for every $n\in\N$, $t\in[0,\thi]$ and compact set $L\subseteq\C$.

So, for specific compact set $L\subseteq\C$ we have $G(z)=G(z_0)$ for every $z\in L$ thus function $G$ is a constant function $G(z)=G(z_0)$, for every $z\in\C$, which is false because $G'(0)\neq0$, according to our hypothesis.

So, we have $a_n\neq0$ for every $n\in\N$. We suppose that the set $\{m_n,n\in\N\}$ is finite. Then, we have that the set $\{a_n,n\in\N\}$ is finite and because $a_n\ra0$ we get that there is some $\nu_0\in\N$, so that $a_n=a_{\nu_0}=0$ for every $n\in\N$, $n\ge \nu_0$, that is false. Thus, the set $\{m_n,n\in\N\}$ is infinite and this implies that there is a subset $\{m'_n,n\in\N\}\subseteq\{m_n,n\in\N\}$ , so that $m'_n$, $n\in\N$ to be a sequence of $\la$, from different terms. \qb
\end{Proof}
\begin{lem}\label{lem3.2}
Let $m=(m_n)_{n\in\N}$ be a sequence of complex numbers, $\thi$ be a positive number and $f,g$ be two entire functions. We suppose that for every compact subset $L$ of $\C$ we have:
\[
\sup_{(t,z)\in[0,\thi]\times L}|f(z+m_ne^{2\pi it})-g(z)|\ra0 \ \ \text{as} \ \ n\ra+\infty.
\]
Then, for every compact subset $L$ of $\C$:
\[
\sup_{(t,z)\in[0,\thi]\times L}|f'(z+m_ne^{2\pi it})-g'(z)|\ra0 \ \ \text{as} \ \ n\ra\infty.
\]
\end{lem}
\begin{Proof}
We fix some compact subset $L$ of $\C$. Let $n_0\in\N$, so that $L$ is a subset of $\oD(0,n_0)$, where $D(0,n)=\{z\in\C|\,|z|<n\}$ for every $n\in\N$.

It is easy to see that
\setcounter{equation}{0}
\begin{eqnarray}
\bigcup_{z\in\oD(0,n_0)}\oD\bigg(z,\frac{n_0}{3}\bigg)\subseteq\oD(0,2n_0).  \label{eq1}
\end{eqnarray}
We shall consider the sequence of functions $F_n$, $n=1,2,\ld$, $F_n:[0,\thi]\times\C\ra\C$, so that $F_n((t,z))=f(z+m_ne^{2\pi it})-g(z)$, for every $(t,z)\in[0,\thi]\times\C$, $n\in\N$ and their partial functions
\[
F^t_n(z)=F_n((t,z)) \ \ \text{for every} \ \ (t,z)\in[0,\thi]\times\C.
\]
Based on the hypothesis, we conclude that for every compact subset $K\subseteq\C$ we have:
\[
\parallel F_n\parallel_{[0,\thi]\times K}\ra0 \ \ \text{as} \ \ n\ra+\infty
\]
where
\[
\parallel F_n\parallel_{[0,\thi]\times K}=\sup_{w\in[0,\thi]\times K}|F_n(w)| \ \ \text{for} \ \ n=1,2,\ld\;.
\]
Let some $\e_0>0$. We set $\e_1=\dfrac{n_0\e_0}{3}$.

Based on our hypothesis, there is a natural number $\nu_0\in\N$, so that for every $n\in\N$, $n\ge\nu_0$, $\parallel F_n\parallel_{[0,\thi]\times\oD_{2n_0}}<\e_1$, where we applied our hypothesis for
\begin{eqnarray}
K=\oD_{2n_0}.  \label{eq2}
\end{eqnarray}
We fix some $t_0\in[0,\thi]$. Then, function $F^{t_0}_n$ is entire for every $n\in\N$.

Let some $z\in\oD_{n_0}$.

Based on Canchy's estimates we have:
\begin{eqnarray}
|(F^{t_0}_n)'(z)|\le\frac{3}{n_0}\sup_{w\in\oD\big(z,\frac{n_0}{3}\big)}|F^{t_0}_n(w)|\le
\frac{3}{n_0}\parallel F_n\parallel_{[0,\thi]\times\oD_{2n_0}}, \label{eq3}
\end{eqnarray}
where for the second inequality we used relation (\ref{eq1}).

Based on inequality (\ref{eq3}), we have:
\begin{eqnarray}
\sup_{(t,z)\in[0,\thi]\times L}|f'(z+m_ne^{2\pi it})-g'(z)|\le\frac{3}{n_0}\parallel F_n\parallel_{[0,\thi]\times\oD_{2n_0}}, \ \ \text{for every} \ \ n\in\N.  \label{eq4}
\end{eqnarray}
Based on (\ref{eq2}) and (\ref{eq4}) we have that for every $n\in\N$, $n\ge\nu_0$:
\begin{eqnarray}
\sup_{(t,z)\in[0,\thi]\times L}|f'(z+m_ne^{2\pi it})-g'(z)|<\e_0.  \label{eq5}
\end{eqnarray}
This gives that:
\[
\sup_{(t,z)\in[0,\thi]\times L}|f'(z+m_ne^{2\pi it})-g'(z)|\ra0 \ \ \text{as} \ \ n\ra+\infty,
\]
that implies the desired result for every compact subset $L$ of $\C$. \qb
\end{Proof}
\begin{lem}\label{lem3.3}
Let $K$ be a compact subset of $\C$ and $f:K\times\C\ra\C$ be a continuous function. Let $f_t:\C\ra\C$ be the partial function of $f$ for every $t\in K$, so that $f_t(z)=f((t,z))$ for every $z\in\C$. We suppose that the partial functions $f_t$ are entire for every $t\in K$.

We shall assume function: $\vPsi:K\times\C\ra\C$ so that:
\[
\vPsi((t,z))=\int_{[0,z]}f_t(j)dj \ \ \text{for every} \ \ (t,z)\in K\times\C.
\]
Then, function $\vPsi$ is continuous.
\end{lem}
\begin{Proof}
We suppose that the compact set is endowed with the relative metric of the usual metric of $\C$ and $K\times\C$ is endowed with the produced topology of the previous topology of $K$ and the usual topology of $\C$.

We set $w_0=(t_0,z_0)\in K\times\C$.

We shall prove that $\vPsi$ is continuous in $w_0$. Let $w_n=(t_n,z_n)\in K\times\C$ for $n=1,2,\ld$ so that $w_n\ra w_0$. It suffices to prove that $\vPsi(w_n)\ra\vPsi(w_0)$.

Given that $w_n\ra w_0$ we have $t_n\ra t_0$ and $z_n\ra z_0$ as $n\ra+\infty$. Given that $z_n\ra z_0$, the sequence $(z_n)_{n\in\N}$ is bounded, so there is a $\nu_0\in\N$, so that
\[
\{z_n,n=1,2,\ld\}\cup\{z_0\}\subseteq\oD_{\nu_0}=\overline{D(0,\nu_0)}=\{z\in\C|\,|z|\le\nu_0\}.
\]
The set $K\times\oD_{\nu_0}$ is a compact subset of $K\times\C$ and function $f$ is continuous on $K\times\oD_{\nu_0}$ and consequently there is some positive number $M>0$, so that:
\setcounter{equation}{0}
\begin{eqnarray}
\sup_{w\in K\times\oD_{\nu_0}}|f(w)|=\parallel f\parallel_ {K\times\oD_{\nu_0}}<M.  \label{eq1}
\end{eqnarray}
For every $n\in\N$ we have:
\begin{align}
|\vPsi(w_n)-\vPsi(w_0)|&=|\vPsi((t_n,z_n))+\vPsi((t_0,z_0))| \nonumber \\
&=|\vPsi((t_n,z_n))-\vPsi((t_n,z_0))+\vPsi((t_n,z_0))-\vPsi((t_0,z_0))| \nonumber \\
&\le|\vPsi((t_n,z_n))-\vPsi((t_n,z_0))|+|\vPsi((t_n,z_0))-\vPsi((t_0,z_0))|.  \label{eq2}
\end{align}
We shall consider the oriented segments $[0,z_0]$, $[z_0,z_n]$ and $[z_n,0]$ of $\C$ and the curve $\ga=[0,z_0]\cup[z_0,z_n]\cup[z_n,0]$ of $\C$ (of course, some of them can have only a single point).

Based on the properties of complex integrals, we have:
\begin{eqnarray}
\int_\ga f_{t_n}(j)dj=\int_{[0,z_0]}f_{t_n}(j)dj+\int_{[z_0,z_n]}f_{t_n}(j)dj+\int_{[z_n,0]}f_{t_n}(j)dj. \label{eq3}
\end{eqnarray}
If the points $0,z_0,z_n$ are not in the same line, we have:
\begin{eqnarray}
\int_\ga f_{t_n}(j)dj=0.  \label{eq4}
\end{eqnarray}
Equality (\ref{eq4}) is supported by Cauchy's Theorem for a Triangle (see Rudin), whereas if the points $0,z_0,z_n$ are in the same line, we can also use equality (\ref{eq4}), according to the properties of complex integrals. Based on equalities (\ref{eq3}) and (\ref{eq4}), we have:
\begin{eqnarray}
\int_{[0,z_n]}f_{t_n}(j)dj-\int_{[0,z_0]}f_{t_n}(j)dj=\int_{[z_0,z_n]}f_{t_n}(j)dj.  \label{eq5}
\end{eqnarray}
Based on equality (\ref{eq5}) we have:
\begin{eqnarray}
\bigg|\int_{[0,z_n]}f_{t_n}(j)dj-\int_{[0,z_0]}f_{t_n}(j)dj\bigg|\le|z_n-z_0|\sup_{z\in[z_0,z_n]}|f_{t_n}(z)|. \label{eq6}
\end{eqnarray}
According to inequality (\ref{eq6}), our hypothesis that $\{z_n,n=1,2,\ld\}\cup\{z_0\}\subseteq\oD_{\nu_0}$ and inequality (\ref{eq1}) we take:
\begin{eqnarray}
\bigg|\int_{[0,z_n]}f_{t_n}(j)dj-\int_{[0,z_0]}f_{t_n}(j)dj\bigg|\le|z_n-z_0|
\cdot M \ \ \text{for every} \ \ n\in\N.  \label{eq7}
\end{eqnarray}
Now for every $n\in\N$ we have
\begin{align}
|\vPsi((t_n,z_0))-\vPsi((t_0,z_0))|&=\bigg|\int_{[0,z_0]}f_{t_n}(j)dj-\int_{[0,z_0]}f_{t_0}(j)dj\bigg| \nonumber\\
&=\bigg|\int_{[0,z_0]}(f_{t_n}(j)-f_{t_0}(j))dj\bigg|\le|z_0|\sup_{j\in[0,z_0]}|f_{t_n}(j)-f_{t_0}(j)| \nonumber \\
&\le\nu_0\parallel f_{t_n}-f_{t_0}\parallel_{\oD_{\nu_0}},  \label{eq8}
\end{align}
because $z_0\in\oD_{\nu_0}$.

We shall prove now that
\[
\para f_{t_n}-f_{t_0}\para_{\oD_{\nu_0}}\ra0 \ \ \text{as} \ \ n\ra+\infty.
\]
Based on our hypothesis, function $f$ is continuous on the compact set $K\times\oD_{\nu_0}$ and thus, it is uniformly continuous.

We fix $\e_0>0$. Provided that $f$ is uniformly continuous on $K\times\oD_{\nu_0}$, there is some $\de_0>0$ so that for every $w_1,w_2\in K\times\oD_{\nu_0}$ with $\para w_1-w_2\para<\de_0$ holds
\begin{eqnarray}
|f(w_1)-f(w_2)|<\e_0 \label{eq9}
\end{eqnarray}
where for every $a=(b,c)\in\C\times\C$, $\para a\para=\sqrt{|b|^2+|c|^2}$, and if $c=(c_1,c_2)$, $b=(b_1,b_2)$ then $|c|=\sqrt{c^2_1+c^2_2}$,
$|b|=\sqrt{b^2_1+b^2_2}$, $a\in\C\times\C$.

Given that $t_n\ra t_0$ based on our hypothesis, there is some $n_0\in\N$, so that for every $n\in\N$, $n\ge n_0$
\begin{eqnarray}
|t_n-t_0|<\de_0.  \label{eq10}
\end{eqnarray}
Based on (\ref{eq9}), (\ref{eq10}) and the fact that functions $f_{t_n}-f_{t_0}$ are continuous on $\C$ we take
\[
\para f_{t_n}-f_{t_0}\para_{\oD_{\nu_0}}<\e_0 \ \ \text{for every} \ \ n\in\N, \ \ n\ge n_0.
\]
As a result
\begin{eqnarray}
\para f_{t_n}-f_{t_0}\para_{\oD_{\nu_0}}\ra0 \ \ \text{as} \ \ n\ra+\infty.  \label{eq11}
\end{eqnarray}
Now, based on (\ref{eq2}), (\ref{eq7}), (\ref{eq8}) and (\ref{eq11}) we take that $\vPsi(w_n)\ra\vPsi(w_0)$ and the proof of this lemma is complete. \qb
\end{Proof}

Next, we use some notations. Let $g\in\ch(\C)$ be an entire function. We assume that $G:\C\ra\C$, where $G(z)=\int_{[0,z]}g(j)dj$, for every $z\in\C$.

Of course, $G$ is an entire function as it is already known by complex analysis.

Let $K$ be a non-empty compact subset of $\C$ and $\Fi_n:K\times\C\ra\C$ be a sequence of continuous functions for $n=1,2,\ld$ so that for every $a\in K$ the functions $\Fi^a_n:\C\ra\C$ where $\Fi^a_n(z)=\Fi_n((a,z))$ for every $z\in\C$, are entire. We also have the functions: $\vPsi_n:K\times\C\ra\C$, so that $\vPsi_n((a,z))=\int_{[0,z]}\Fi^a_n(j)dj$, for every $n\in\N$, $(a,z)\in K\times\C$.

Functions $\Fi^a_n$, $n\in\N$, $a\in K$ are partial functions of $\Fi_n$, $n\in\N$.

Functions $\vPsi_n$, $n\in\N$ are continuous by Lemma \ref{lem3.3}. We shall consider the partial functions $\vPsi^a_n:\C\ra\C$, so that $\vPsi^a_n(z)=\vPsi_n((a,z))$, for every $n\in\N$, $a\in K$ and $z\in\C$. It is well-known by complex analysis that functions $\vPsi^a_n$ are entire for every $a\in K$ and $n\in\N$.

Upon defining the above notations, we shall continue with the following Lemma \ref{lem3.4}.
\begin{lem}\label{lem3.4}
We suppose that for every compact subset $L\subseteq\C$ the following shall apply: $(\ast)$
\[
\sup_{(a,z)\in K\times L}|\Fi_n((a,z))-g(z)|\ra0 \ \ \text{as} \ \ n\ra+\infty.
\]
Then, for every compact subset $L\subseteq\C$, it shall apply:
\[
\sup_{(a,z)\in K\times L}|\vPsi_n((a,z))-G(z)|\ra0 \ \ \text{as} \ \ n\ra+\infty.
\]
\end{lem}
\begin{Proof}
Let $\Ga:=\{\vPsi^a_n|a\in K,n\in\N\}$.

It is obvious that $\Ga\subseteq\ch(\C)$.

We shall firstly prove that $\Ga$ is a relatively compact subset of $\ch(\C)$, endowed with the topology of uniform convergence on compacta.

We write $D_n=\{z\in\C|\,|z|<n\}$, for every $n\in\N$.

We fix $n_0\in\N$. Based on our hypothesis $(\ast)$, there is some natural number $m_0(n_0)\in\N$, $m_0(n_0)\ge2$, so that for every $n\in\N$, $n\ge m_0(n_0)$ the following shall apply:
\setcounter{equation}{0}
\begin{eqnarray}
\sup_{(a,z)\in K\times\oD_{n_0}}|\Fi_n((a,z))-1|<1.  \label{eq1}
\end{eqnarray}
Based on (\ref{eq1}), we take that for every $n\in\N$, $n\ge m_0(n_0)$,
\begin{eqnarray}
\para\Fi_n\para_{K\times\oD_{n_0}}<2    \label{eq2}
\end{eqnarray}
where
\[
\para\Fi_n\para_{K\times\oD_{n_0}}=\sup_{(a,z)\in K\times\oD_{n_0}}|\Fi_n((a,z))|, \ \ \text{for every} \ \ n\in\N.
\]
Let
\[
M_1:=\max\{\para\Fi_n\para_{K\times\oD_{n_0}}, \; n=1,\ld,m_0(n_0)-1\}.
\]
We set $M_2=\max\{M_1,2\}+1$.

Then, based on (\ref{eq2}) and the definitions of numbers $M_1,M_2$, we have:
\begin{eqnarray}
\para\Fi_n\para_{K\times\oD_{n_0}}<M_2, \ \ \text{for every} \ \ n=1,2,\ld\;.  \label{eq3}
\end{eqnarray}
According to the definition of functions $\vPsi^\al_n$, $\al\in K$, $n\in\N$ and inequality (\ref{eq3}) we have:
\begin{eqnarray}
\para\vPsi^\al_n\para_{\oD_{n_0}}<n_0\cdot M_2 \ \ \text{for every} \ \ \al\in K \ \ \text{and} \ \ n\in\N.  \label{eq4}
\end{eqnarray}
Inequality (\ref{eq4}) gives us that for every compact subset $L\subseteq\C$, there is some constant $M_L>0$, so that
\begin{eqnarray}
\text{For every} \ \ \al\in K, \; n\in\N, \; z\in L, \ \ |\vPsi^\al_n(z)|<M_L.  \label{eq5}
\end{eqnarray}
Inequality (\ref{eq5}) gives us that set $\Ga$ is a locally bounded subset of $\ch(\C)$. By Montel's Theorem we take that set $\Ga$ is a relatively compact subset of $\ch(\C)$ and so, based on Arzela-Ascoli Theorem, we take that set $\Ga$ is equicontinuous. We set $\mathfrak{D}:=\{\vPsi_n,n\in\N\}$. \\
We shall prove now that set $\mfd$ is equicontinuous.

We fix $w_0\in K\times\C$.

We shall prove that for every $\e>0$ there is $\de>0$ so that for every $w\in K\times\C$, $\rho((w,w_0))<\de\; (w\in N_{w_0}(\de))$ for every $n\in\N$
\begin{eqnarray}
|\vPsi_n(w)-\vPsi_n(w_0)|<\e.  \label{eq6}
\end{eqnarray}
We fix $\e_0>0$. Let $w_0=(a_0,z_0)$, where $a_0\in K$ and $z_0\in\C$.

We fix a $\nu_0\in\N$, so that: $\nu_0>|z_0|+1$.

Based on our hypothesis $(\ast)$ there is some $m_1\in\N$, $m_1\ge2$, so that for every $n\in\N$, $n\ge M_1$
\begin{eqnarray}
\sup_{(\al,z)\in K\times\oD_{\nu_0}}|\Fi_n((a,z))-g(z)|<\frac{\e_0}{4\nu_0}.  \label{eq7}
\end{eqnarray}
Now, based on (\ref{eq7}) and the definition of functions $\vPsi_n$, $n\in\N$, $n\ge m_1$, we take that for every $\al\in K$, $z\in\oD_{\nu_0}$ the following applies:
\begin{align}
|\vPsi_n((\al,z))-\vPsi_n(\al_0,z)|&\le\nu_0\cdot\sup_{z\in\oD_{\nu_0}}|\Fi_n((\al,z))-\Fi_n((\al_0,z)) \nonumber \\
&\le2\nu_0\sup_{(\al,z)\in K\times\oD_{\nu_0}}|\Fi_n((\al,z))-g(z)|<\frac{\e_0}{2}.  \label{eq8}
\end{align}
Based on the equicontinuity of $\Ga$ on $z_0$ we take that for $\e_0/4$ there is $\de_1\in(0,1)$ so that for every $z\in\C$, $|z-z_0|<\de_1$,
\begin{eqnarray}
|\vPsi^\al_n(z)-\vPsi^\al_n(z_0)|<\frac{\e_0}{2} \ \ \text{for every} \ \ \al\in K \ \ \text{and} \ \ n=1,2,\ld\;.  \label{eq9}
\end{eqnarray}
It is obvious that for every $z\in\C$, $|z-z_0|<\de_1$ we have
\begin{eqnarray}
z\in\oD_{\nu_0}.  \label{eq10}
\end{eqnarray}
According to (\ref{eq10}) and (\ref{eq8}), we conclude that for every $z\in\C$, $|z-z_0|<\de_1$ for every $a\in K$ and $n\in\N$, $n\ge m_1$ the following applies:
\begin{eqnarray}
|\vPsi_n((a,z))-\vPsi_n((a_0,z))|<\frac{\e_0}{2}.  \label{eq11}
\end{eqnarray}
According to (\ref{eq9}) and (\ref{eq11}) we take that for every $a\in K$,. $z\in\C$, $|z-z_0|<\de_1$ and $n\in\N$, $n\ge m_1$, we have:
\begin{eqnarray}
|\vPsi_n((a,z))-\vPsi_n((a_0,z_0))|<\e_0.  \label{eq12}
\end{eqnarray}
Based on (\ref{eq12}) and the fact that functions $\vPsi_n$ are continuous on $w_0$ for $n=1,\ld,m_1-1$, we conclude that there is some $\de_0\le\de_1$, so that for every $w\in K\times\C$ and $\rho((w,w_0))<\de_0$ the following applies:
\begin{eqnarray}
|\vPsi_n(w)-\vPsi_n(w)|<\e_0, \ \ \text{for every} \ \  n=1,2,\ld,  \label{eq13}
\end{eqnarray}
where $\rho:(K\times\C)\times(K\times\C)\ra\R$ is the usual $\rho$ metric on $K\times\C$, where if $(a_1,z_1),(a_2,z_2)\in K\times\C$
\[
\rho((a_1,z_1),(a_2,z_2))=\sqrt{|a_2-a_1|^2+|z_1-z_2|^2},
\]
condition (\ref{eq13}) tells us that set $\mfd$ is equicontinuous on fixed $w_0\in K\times\C$, and so set $\mfd$ is equicontinuous on $K\times\C$.

We shall prove now that for every $(a,z)\in K\times\C$ holds
\[
\vPsi_n((a,z))\ra G(z) \ \ \text{as} \ \ n\ra+\infty.
\]
We fix $w_1=(a_1,z_1)\in K\times\C$.

Let $n_1\in\N$, so that $z_1\in\oD_{n_1}$.

According to our hypothesis $(\ast)$, we have:
\[
\sup_{(a,z)\in K\times\oD_{n_1}}|\Fi_n((a,z))-g(z)|\ra0 \ \ \text{as} \ \ n\ra+\infty,
\]
and this gives that $\Fi^a_n\ra g$ uniformly on $[0,z_1]$ because
\begin{eqnarray}
z_1\in\oD_{n_1},  \label{eq14}
\end{eqnarray}
By convergence of (\ref{eq14}) we have:\\
$\vPsi^a_n(z_1)\ra G(z_1)$ and this gives that
\begin{eqnarray}
\vPsi_n((a,z))\ra G(z), \ \ \text{for every} \ \ (a,z)\in K\times\C.  \label{eq15}
\end{eqnarray}
We shall consider the functions $\vPsi_{m,n}:K\times\oD_m\ra\C$, so that
\[
\vPsi_{m,n}((a,z))=\vPsi_n((a,z)) \ \ \text{for every} \ \ n\in\N, \ \ (a,z)\in K\times\oD_m,
\]
that is $\vPsi_{m,n}=\vPsi\upharpoonright K\times\oD_m$ is the restriction of $\vPsi_n$ on $K\times\oD_m$ for every $m\in\N$.

We set $\mfd_m=\{\vPsi_{m,n}|n\in\N\}$, for every $m\in\N$. Because $\mfd$ is equicontinuous on $K\times\C$ we take that the sets $\mfd_m$ are equicontinuous on $K\times\oD_m$ for every $m\in\N$.

Given that $\Ga$ is locally bounded we take that for every $m\in\N$, there is $M_m>0$ so that:
\[
\para\vPsi_n\para_{K\times\oD_m}<M_m, \ \ \text{for every} \ \ n=1,2,\ld\;.
\]
This gives that $\vPsi_{m,n}\in B(0,M_n)$ for every $n=1,2,\ld$, on $C(K\times\oD_m)$ for every $m\in\N$.

So, we have $\mfd_m\subseteq B(0,M_m)$ for every $m\in\N$. Therefore, the $\mfd_m$ is a bounded subset on metric space $(C(K\times\oD_m,\para\cdot\para\infty))$. So, the $\overline{\mfd}_m$ set is a closed and bounded subset of $C(K\times\oD_m),\para\cdot\para\infty))$, for every $m\in\N$).

We shall prove now that the set $\overline{\mfd}_m$ is equicontinuous on $K\times\oD_m$, for every $m\in\N$.

We fix $m_2\in\N$ and $w_2\in K\times\oD_{m_2}$. We fix $\e_1>0$. We shall prove that $\overline{\mfd}_{m_2}$ is equicontinuous on $w_2$. Because of the equicontinuity of $\mfd_{m_2}$ on $w_2$ we take that for $\frac{\e_1}{3}>0$ there is some $\de_2>0$, so that for every $w\in K\times\oD_{m_2}$, $\rho((w,w_2))<\de_2$ we have:
\begin{eqnarray}
|\vPsi_{m_2,n}(w)-\vPsi_{m_2,n}(w_0)|<\frac{\e_1}{3} \ \ \text{for every} \ \ n=1,2,\ld\;.  \label{eq16}
\end{eqnarray}
Let $h\in\overline{\mfd}_{m_2}$. Then, there is subsequence $(k_n)_{n\in\N}$ of natural numbers, so that $\vPsi_{m_2,k_n}\ra h$ uniformly on $K\times\oD_{m_2}$ as $n\ra+\infty$.

Thus, there is $m_3\in\N$, so that for every $n\in\N$, $n\ge m_3$ has as follows:
\begin{eqnarray}
\para\vPsi_{m_2,k_2}-h\para_{k\times\oD_{m_2}}<\frac{\e_1}{3}.  \label{eq17}
\end{eqnarray}
We apply now (\ref{eq16}) for $n=k_{m_3}$ and (\ref{eq17}) for $n=m_3$ we take that for every $w\in K\times\oD_{m_2}$, $\rho((w,w_2))<\de_2$ we have:
\[
|h(w)-h(w_2)|=|h(w)-\vPsi_{m_2,k_{m_3}}(w)+\vPsi_{m_2,k_{m_3}}(w)-\vPsi_{m_2,k_{m_3}}(w_2)+
\vPsi_{m_2,k_{m_3}}(w_2)-h(w_2)|<\e_1.
\]
Because of the fact that number $\de_2$ depends only on $m_2$, $w_2$, $\e_1$ and not on $h$, we take that the set $\overline{\mfd}_{m_2}$ is equicontinuous on $w_2$, so the set $\overline{\mfd}_{m_2}$ is equicontinuous on $K\times\oD_{m_2}$.

We have proved that the set $\overline{\mfd}_{m_2}$ is a closed, bounded and equicontinuous subset of $(C(K\times\oD_{m_2}),\para\cdot\para\infty)$. Based on Ascoli's-Arzela Theorem, we take that the set $\overline{\mfd}_{m_2}$ is a compact subset of $C(K\times\oD_{m_2})$. Thus, we take that the set $\overline{\mfd}_m$ is a compact subset of $C(K\times\oD_m)$ for every $m\in\N$.

We set $\oG:K\times\C\ra\C$, where $\oG((a,z))=G(z)$ for every $(a,z)\in K\times\C$. Of course, $\oG\in C(K\times\C)$.

We shall prove now that $\vPsi_n\ra\oG$ on the space $C(K\times\C)$ endowed with the topology of uniform convergence on compact subsets of $K\times\C$.

For the purposes of proof by contradiction, we suppose that the sequence $(\vPsi_n)_{n\in\N}$ does not converge to $\oG$ on the space $C(K\times\C)$.

Then, there is some $\nu_1\in\N$, some positive number $\e_2>0$ and some subsequence $(a_n)_{n\in\N}$ of natural numbers, so that:
\begin{eqnarray}
\para\vPsi_{a_n}-\oG\para_{K\times\oD_{\nu_1}}\ge\e_2 \ \ \text{for every} \ \ n=1,2,\ld\;.  \label{eq18}
\end{eqnarray}
The sequence $(\vPsi_{\nu_1,a_n})_{n\in\N}$ is a sequence in $\overline{\mfd}_{\nu_1}$. Because of the fact that $\overline{\mfd}_{\nu_1}$ is compact, sequence $(\vPsi_{\nu_1,a_n})_{n\in\N}$ has a subsequence $(\vPsi_{\nu_1,a_{\rho_n}})_{n\in\N}$ that converges to some\linebreak $h_1\in C(K\times\oD_{\nu_1})$.

Therefore, we have:
\begin{eqnarray}
\vPsi_{\nu_1,a_{\rho_n}}\ra h_1 \ \ \text{uniformly on} \ \ K\times\oD_{\nu_1}.  \label{eq19}
\end{eqnarray}
Based on (\ref{eq15}), we have
\[
\vPsi_n((a,z))\ra G(z) \ \ \text{for every} \ \ (a,z)\in K\times\C.
\]
This gives us:
\begin{eqnarray}
\vPsi_{a_{\rho_n}}(w)\ra\oG(w) \ \ \text{for every} \
 \ w\in K\times\oD_{\nu_1}.  \label{eq20}
\end{eqnarray}
Based on (\ref{eq19}), we have:
\begin{eqnarray}
\vPsi_{\nu_1,a_{\rho_n}}(w)\ra h_1(w) \ \ \text{for every} \ \ w\in K\times\oD_{\nu_1}.  \label{eq21}
\end{eqnarray}
Based on (\ref{eq20}) and (\ref{eq21}), we have:
\begin{eqnarray}
h_1(w)=\oG(w) \ \ \text{for every} \ \ w\in K\times\oD_{\nu_1}.  \label{eq22}
\end{eqnarray}
Based on (\ref{eq19}) and (\ref{eq22}), we have:
\begin{eqnarray}
\vPsi_{\nu_1,a_{\rho_n}}\ra\oG \ \ \text{uniformly on} \ \ K\times\oD_{\nu_1}. n \label{eq23}
\end{eqnarray}
Based on (\ref{eq18}) and (\ref{eq23}) we end up to a contradiction. So, we have proved that $\vPsi_n\ra\oG$ uniformly on compact subsets of $K\times\C$, or else:
\[
\sup_{(a,z)\in K\times L}|\vPsi_n((a,z))-G(z)|\ra0 \ \ \text{as} \ \ n\ra+\infty
\]
for very compact set $L\subseteq\C$. \qb
\end{Proof}
\begin{lem}\label{lem3.5}
Let $(m_n)_{n\in\N}$ be a sequence of complex numbers, so that $m_n\ra\infty$, $\thi\in(0,1)$ and $a\in\C$, $a\neq0$.

Then, there is no entire function $f$, so that:
\[
\sup_{t\in[0,\thi]}|f(m_ne^{2\pi it})-a|\ra0 \ \ \text{and}
\]
\[
\int_{[m_n,m_ne^{2\pi i\thi}]}f(z)dz\ra0 \ \ \text{as} \ \ n\ra+\infty.
\]
\end{lem}
\begin{Proof}
To take a contradiction we suppose that there exists an entire function $f$ that satisfies the above two convergence.

We have the curves $\ga_n:[0,\thi]\ra\C$, where $\ga_n(t)=m_ne^{2\pi it}$ for every $t\in[0,\thi]$, $n\in\N$. We also have $\ga^\ast_n:=\ga_n([0,\thi])$ for $n=1,2,\ld\;.$

Because $m_n\ra\infty$, we use only the terms $m_n$, $n\in\N$, such that $m_n\neq0$.

Based on Cauchy's Theorem we have:
\setcounter{equation}{0}
\begin{eqnarray}
\int_{[m_n,m_ne^{2\pi i\thi}]}f(z)dz=\int_{\ga^\ast_n}f(z)d(z)  \ \ \text{for every} \ \ n=1,2,\ld\;. \label{eq1}
\end{eqnarray}
We have: $A:\C\ra\C$ for the constant function, so that $A(z)=a$ for every $z\in\C$.

We also have:
\begin{eqnarray}
\int_{[m_n,m_ne^{2\pi i\thi]}}A(z)dz=a(m_ne^{2\pi i\thi}-m_n) \ \ \text{for every} \ \ n=1,2,\ld\;.  \label{eq2}
\end{eqnarray}
We fix $\e_0\in\Big(0,\dfrac{|a|\cdot|e^{2\pi i\thi}-1|}{2\pi}\Big)$.

We can write down the first of the two convergences of hypothesis as follows:
\begin{eqnarray}
\sup_{z\in\ga^\ast_n}|f(z)-a|\ra0 \ \ \text{as} \ \ n\ra+\infty.  \label{eq3}
\end{eqnarray}
By our hypothesis (\ref{eq1}) and (\ref{eq3}) we take that there is some $\nu_0\in\N$, so that for every $n\in\N$, $n\ge\nu_0$ has as follows:
\begin{eqnarray}
\sup_{z\in\ga^\ast_n}|f(z)-a|<\e_0 \ \ \text{and}  \label{eq4}
\end{eqnarray}
\begin{eqnarray}
\bigg|\int_{\ga^\ast_n}f(z)dz\bigg|<\e_0, \ \ m_n\neq0.  \label{eq5}
\end{eqnarray}
Based on (\ref{eq2}), Cauchy's Theorem, (\ref{eq4}) and the simple properties of the complex integrals, we have:
\begin{eqnarray}
\bigg|\int_{\ga^\ast_n}f(z)dz-a(m_ne^{2\pi i\thi}-m_n)\bigg|<2\pi|m_n|\e_0.  \label{eq6}
\end{eqnarray}
Based on (\ref{eq4}), (\ref{eq5}), (\ref{eq6}), triangle inequality and the specific of $\e_0$ we assume that for every $n\in\N$, $n\ge\nu_0$, the following applies:
\begin{eqnarray}
|m_n|<\frac{\e_0}{|a|\,|e^{2\pi i\thi}-1|-2\pi\e_0}  \label{eq7}
\end{eqnarray}
(where $|a|\,|e^{2\pi i\thi}-1|-2\pi\e_0>0$ from the certain choice of $\e_0$).

Inequality (\ref{eq7}) and the fact that $m_n\ra\infty$ gives a contradiction and this completes the proof of this lemma. \qb
\end{Proof}
\section{The main result}\label{sec4}
\noindent

In order to prove the main result, we also need the following corollary of Lemma \ref{lem3.2}.
\begin{cor}\label{cor4.1}
Let $(m_n)_{n\in\N}$ be a sequence of complex numbers, $\thi$ be a positive number and $f$, $g$ be two entire functions.

We suppose that for every compact subset $L$ of $\C$ the following shall apply:
\[
\sup_{(t,z)\in[0,\thi]\times L}|f(z+m_ne^{2\pi it})-g(z)|\ra0 \ \ \text{as} \ \ n\ra+\infty.
\]
Then, for every compact subset $L$ of $\C$ and $\nu\in\N$
\[
\sup_{(t,z)\in[0,\thi]\times L}|f^{(\nu)}(z+m_ne^{2\pi it})-g^{(\nu)}(z)|\ra0 \ \ \text{as} \ \ n\ra+\infty.
\]
\end{cor}
\begin{Proof}
It is simple implication of Lemma \ref{lem3.2} by induction.

Now, we are ready to prove the main result of this article. \qb
\end{Proof}
\begin{thm}\label{thm4.2}
Let $\la=(\la_n)_{n\in\N}$ be a sequence of non-zero complex numbers, so that $\la_n\ra\infty$, $\thi\in(0,1)$ and $G\in\ch(\C)$, where $G$ is not a constant function.

Then, we have: $\cu(\la,\thi,G)=\emptyset$.
\end{thm}
\begin{Proof}
We shall prove the Theorem by distinguishing some cases.
\begin{itemize}
        \item {\bf Case 1} \\
        $G(0)=0$ and $G'(0)\neq0$.\\
        The result is supported by Proposition 1.
        \item {\bf Case 2} \\
         $G(0)=0$ and $G'(0)=0$.
      \end{itemize}
We shall distinguish two cases here: \smallskip\\
a) $G^{(\nu)}(0)=0$ for every $\nu\in\N$.\\
Provided that $G\in\ch(\C)$ we have $G(z)=\dis\sum^{+\infty}_{\nu=0}\dfrac{G^{(\nu)}(0)}{\nu!}z^\nu$ for every $z\in\C$, so we have $G(z)=0$ for every $z\in\C$, which is false because $G$ is not a constant function in our hypothesis. \smallskip \\
b) There is a $\nu\in\N$, $\nu\ge2$ so that $G^{(\nu)}(0)\neq0$.\\
Let $\nu_0=\min\{\nu\in\N|G^{(\nu)}(0)\neq0\}$, that is $\nu_0$ is the smallest natural number, so that $G^{(\nu_0)}(0)\neq0$. Of course, $\nu_0\ge2$.

Therefore, we have $G^{(\nu_0-1)}(0)=0$ and $G^{(\nu_0)}(0)\neq0$.

We suppose that $\cu(\la,\thi,G)\neq\emptyset$. Let $f\in\cu(\la,\thi,G)$. Then, there is a sequence $(m_n)_{n\in\N}$, so that $m_n\in\{\la_n,n\in\N\}$ for every $n\in\N$, where for every compact subset $L\subseteq\C$
\[
\sup_{(t,z)\in[0,\thi]\times L}|f(z+m_ne^{2\pi it})-G(z)|\ra0 \ \ \text{as} \ \ n\ra+\infty.
\]
Based on the above Corollary \ref{cor4.1}, we take that for every compact subset $L\subseteq\C$
\setcounter{equation}{0}
\begin{eqnarray}
\sup_{(t,z)\in[0,\thi]\times L}|f^{(\nu_0-1)}(z+m_ne^{2\pi it})-G^{(\nu_0-1)}(z)|\ra0 \ \ \text{as} \ \ n\ra+\infty.  \label{eq1}
\end{eqnarray}
Because $G^{(\nu_0-1)}(0)=0$ and $G^{(\nu_0)}(0)\neq0$ we take that the function $G^{(\nu_0-1)}\in\ch(\C)$ is not a constant function. Of course, the function $f^{(\nu_0-1)}\in\ch(\C)$.

Based on (\ref{eq1}) and Proposition \ref{prop2.1} we have a contradiction, because, according to (\ref{eq1}), we have $f^{(\nu_0-1)}\in\cu(\la,\thi,G^{(\nu_0-1)})$, that is $\cu(\la,\thi,G^{(\nu_0-1)})\neq\emptyset$ that is false by Proposition \ref{prop2.1}.
\begin{itemize}
  \item {\bf Case 3} \\
  $G(0)\neq0$
\end{itemize}
We have function $g=G-G(0)$. Of course, we have $g\in\ch(\C)$ and $g(0)=0$.

We suppose that $\cu(\la,\thi,G)\neq\emptyset$. Let $f\in U(\la,\thi,G)$.

Then, there is a sequence $m=(m_n)_{n\in\N}$ of complex numbers, where $m_n\in\{\la_n,n\in\N\}$ for every $n\in\N$, so that for every compact subset $L\subseteq\C$
\begin{eqnarray}
\sup_{(t,z)\in[0,\thi]\times L}|f(z+m_ne^{2\pi it})-G(z)|\ra0 \ \ \text{as} \ \ n\ra+\infty. \label{eq2}
\end{eqnarray}
We have function $F=f-G(0)$. Of course, $F\in\ch(\C)$.

Based on (\ref{eq2}) we have:
\begin{eqnarray}
\sup_{(t,z)\in[0,\thi]\times L}|F(z+m_ne^{2\pi it})-g(z)|\ra0 \ \ \text{as} \ \ n\ra+\infty, \ \ \text{for every compact set} \  L\subseteq\C.  \label{eq3}
\end{eqnarray}
Based on (\ref{eq3}) we have: $F\in\ch(\la,\thi,g)$ where $g$ is not a constant, according to its definition, that is, $\cu(\la,\thi,g)\neq\emptyset$, which is false by the above Cases 1 and 2, because $g(0)=0$.

The proof of our main result is complete now. \qb
\end{Proof}
\vspace*{1cm}
N. Tsirivas  \\
Department of Mathematics, \\
University of Thessaly, \\
Lamia, Greece.\\
email: ntsirivas@uth.gr

\end{document}